\documentclass[12pt,a4paper]{article} 
\usepackage{amssymb,amsthm,amsmath}

\def\otb{\otimes_B}
\def\ot4{\otimes_{\A(\S^4)}}

\renewcommand{\bar}[1]{\overline{#1}}

\def\det{\mathrm{det} \mb}
\def\dd{\mathrm{d}} 

\def\ii{\mathrm{i}} 
\def\id{\mathrm{id}}
\def\into{\hookrightarrow}

\def\tr{\mathrm{Tr}}
\def\t{\mathrm{t}}

\def\Aut{\mathrm{Aut}}

\def\End{\mathrm{End}}
\def\Hom{\mathrm{Hom}}
\def\Mat{\mathrm{Mat}}

\def\isom{\simeq}

\def\bC{\mathbb{C}} 

\def\I{\mathbb{I}} 
 
\def\bR{\mathbb{R}}

\def\bT{\mathbb{T}}
\def\bZ{\mathbb{Z}} 

\def\A{\mathcal{A}} 

\def\cB{\mathcal{B}}

\def\cE{\mathcal{E}}
\def\cH{\mathcal{H}}

\def\cS{\mathcal{S}}

\newcommand{\ket}[1]{|#1\rangle}    
\newcommand{\bra}[1]{\langle#1|}    
\newcommand{\braket}[2]{\langle#1|#2\rangle}
\newcommand{\ketbra}[2]{\langle#1|#2\rangle}

\newcommand{\br}[1]{\langle #1 \rangle}
\def\ch{\mathrm{ch}}

\def\coinv{\mathrm{Coinv}}

\def\dix{\mathrm{Tr}_\omega \mb}
\def\ind{\mathrm{Ind}\, }

\def\inv{\mathrm{Inv}}

\def\res{\mathrm{Res}\mb}
\def\resz{\underset{z=0}{\res}}

\def\Ad{\mathrm{Ad}}
\def\M{M_\theta}
\def\n{{(n)}}
\def\R{\bR_\theta}
\def\S{S_\theta}
\def\Sk{S_{\theta'}}
\def\st{\mathrm{st}} 		
\def\Sym{\mathrm{Sym}}
\def\T{\bT_\theta}
\def\Tk{\bT_{\theta'}}
\def\un{\mathrm{un}}

\def\xc{x_{(0)}}				
\def\zc{z_{(0)}} 				

\def\bd{\begin{displaymath} }
\def\ed{\end{displaymath} }
\def\be{\begin{equation}}
\def\ee{\end{equation}}
\def\bea{\begin{eqnarray}}
\def\eea{\end{eqnarray}}
\def\bean{\begin{eqnarray*}}
\def\eean{\end{eqnarray*}}
\def\bmult{\begin{multline}}
\def\emult{\end{multline}}
\def\nn{\nonumber}
\newtheorem{thm}{Theorem}

\newtheorem{lma}[thm]{Lemma} 
\newtheorem{prop}[thm]{Proposition}
\newtheorem{defn}[thm]{Definition}
\newtheorem{ex}[thm]{Example}
\def\mb{\mbox{ }}
\def\rightbox{\begin{flushright}$\Box$\end{flushright} }

\begin{document}
\title{Principal fibrations from noncommutative spheres}
\author{Giovanni Landi$^1$, Walter van Suijlekom$^2$ \\[10mm]$^1$ Dipartimento di Matematica e Informatica, Universit\`a di Trieste\\Via A. Valerio 12/b, I-34127 
Trieste, Italy\\ 
and INFN, Sezione di Napoli, Napoli, Italy\\
\texttt{landi@univ.trieste.it}\\[5mm]
 $^2$ Scuola Internazionale Superiore di Studi Avanzati\\
Via Beirut 2-4, I-34014 Trieste, Italy\\
\texttt{wdvslkom@sissa.it}
}
\date{}
\maketitle
\hyphenation{de-for-ma-tions}
\begin{abstract}
We construct noncommutative principal fibrations $\S^7 \to \S^4$ which are deformations of the 
classical $SU(2)$ Hopf fibration over the four sphere. 
We realize the noncommutative vector bundles associated to the irreducible representations of $SU(2)$ as modules of coequivariant maps and construct corresponding projections.  The index of Dirac operators with coefficients in the associated bundles is computed with the Connes-Moscovici local index formula.
The algebra inclusion  $\A(\S^4) \into \A(\S^7)$ is an example of a not trivial quantum principal bundle.
\end{abstract}
%
%
%
\thispagestyle{empty}
\newpage
\section{Introduction}\label{sect:intro}
The ADHM construction \cite{Ati79,ADHM78} of instantons in Yang-Mills theory has at its heart the theory of connections on principal and associated bundles. A central example is the basic $SU(2)$-instanton on $S^4$ which is described by the well-known Hopf $SU(2)$-principal bundle $S^7 \to S^4$ and connections thereon. 

In this paper, we consider a noncommutative version of this Hopf fibration, in the framework of the isospectral deformations introduced in \cite{CL01}, while trying to understand the structure behind the noncommutative instanton bundle found there. 

In Sect.~\ref{sect:NCS} we will review the construction of $\theta$-deformed spheres where $\theta$ is an anti-symmetric real-valued matrix. Apart from the noncommutative spheres $\S^m$, we also introduce differential calculi $\Omega(\S^m)$ as quotients of the universal differential calculi. On the sphere $\S^m$ one constructs a noncommutative Riemannian spin geometry $(C^\infty(\S^m), D, \cH)$ in which the Dirac operator $D$ is the classical one and $\cH=L^2(S^m,\cS)$ is the usual Hilbert space of spinors. Then the deformations are isospectral, as mentioned. Furthermore, one also constructs a Hodge star operator $\ast_\theta$ acting on the differential calculus $\Omega(\S^m)$ which is most easily defined using the so-called splitting homomorphism \cite{CD02}.

In Sect.~\ref{sect:hopf}, we focus on two noncommutative spheres $\S^4$ and $\Sk^7$ starting from the algebras $\A(\S^4)$ and $\A(\Sk^7)$ of polynomial functions on them. The latter algebra carries an action of the (classical) group $SU(2)$ by automorphisms in such a way that its invariant elements are exactly the polynomials on $\S^4$. The anti-symmetric $2\times 2$ matrix $\theta$ is given by a single real number also denoted by $\theta$.  
On the other hand, the requirements that $SU(2)$ acts by automorphisms and that  $\S^4$ makes the algebra of invariant functions, give the matrix $\theta'$ in terms of $\theta$. This yields a one-parameter family of noncommutative Hopf fibrations.

For each irreducible representation $V^\n:=\Sym^n(\bC^2)$ of $SU(2)$ we construct the noncommutative vector bundles $E^\n$ associated to the fibration $\Sk^7 \to \S^4$. By dualizing the classical construction, these bundles are described by the module of coequivariant maps from $\bC^2$ to $\A(\Sk^7)$. As expected, these modules are finitely generated projective and we construct explicitly the projections $p_\n \in M_{4^n}(\A(\S^4))$ such that these modules are isomorphic to the image of $p_\n$ in $\A(\S^4)^{4^n}$. Then, one defines connections $\nabla=p_\n \dd$ as maps from $\Gamma(\S^4,E^\n)$ to $\Gamma(\S^4,E^\n) \otimes_{\A(\S^4)} \Omega^1(\S^4)$, where $\Omega^*(\S^4)$ is the quotient of the universal differential calculus mentioned above. The corresponding connection one-form $A$ turns out to be valued in a representation of the Lie algebra $su(2)$. 

By using the projection $p_\n$, the Dirac operator with coefficients in the noncommutative vector bundles $E^\n$ is given by $D_{p_\n}:= p_\n D p_\n$. In order to compute its index, we first show that the local index theorem of Connes and Moscovici \cite{CM95} takes a very simple form in the case of isospectral deformations. Indeed, for these deformations and with any projection $e$, one finds,
\be
\ind D_e = \resz z^{-1} \tr \bigg(\gamma \big(e-\frac{1}{2}\big) |D|^{-2z} \bigg)+\sum_{k \geq 1} c_k \resz \tr \bigg( \gamma \big(e-\frac{1}{2} \big) [D,e]^{2k} |D|^{-2(k+z)} \bigg)
\ee
with some proper coefficients $c_k$. When applied to the projections $p_\n$ on $\S^4$, we obtain exactly as in the classical case, 
\be
\ind D_{p_\n} = \frac{1}{6}n(n+1)(n+2).
\ee

Finally, in Sect.~\ref{sect:hopf-galois} we show that the fibration $\Sk^7 \to \S^4$ is a `not-trivial principal bundle with structure group $SU(2)$'. This means that the inclusion $\A(\S^4) \into \A(\Sk^7)$ is a not-cleft Hopf-Galois extension \cite{KT81,Mon93}; in fact, it is a principal extension \cite{BH04}. On this extension, we find an explicit form of the (strong) connection which induces connections on the associated bundles $E^\n$ as maps from $\Gamma(\S^4,E^\n)$ to $\Gamma(\S^4,E^\n) \otimes_{\A(\S^4)} \Omega^1(\A(\S^4))$, where $\Omega^*(\A(\S^4))$ is the universal differential calculus on $\A(\S^4)$. We show that these connections coincide with the Grassmannian connections $\nabla =p_\n \dd$ on the quotient $\Omega(\S^4)$ of the universal differential calculus alluded to before. 

\section{Noncommutative spherical manifolds} \label{sect:NCS}
In this section, we will recall the construction of the noncommutative spheres $\S^n$ as introduced in \cite{CL01} and elaborated in \cite{CD02}. Essentially, these $\theta$-deformations are a natural extension of the noncommutative torus (for a review see \cite{Rie90a}) to (compact) Riemannian manifolds carrying an action of the $n$-torus $\bT^n$. In this paper we will restrict only to the cases of planes and spheres.

For $\lambda^{\mu\nu}=e^{2\pi \ii \theta_{\mu\nu}}$, where $\theta_{\mu\nu}$ is an anti-symmetric real-valued matrix,
the algebra $\A(\R^{2n})$ of polynomial functions on the noncommutative $2n$-plane is defined to be the unital $*$-algebra generated by $2n$ elements $z^\mu, \bar{z}^\mu (\mu=1,\ldots, n)$ with relations
\be
z^\mu z^\nu = \lambda^{\mu\nu} z^\nu z^\mu;\quad \bar{z}^\mu z^\nu = \lambda^{\nu\mu} z^\nu \bar{z}^\mu;\quad \bar{z}^\mu \bar{z}^\nu = \lambda^{\mu\nu} \bar{z}^\nu \bar{z}^\mu,
\ee
The involution $*$ is defined by putting $z^{\mu*} = \bar{z}^\mu$. For $\theta=0$ one recovers the commutative $*$-algebra of complex polynomial functions on $\bR^{2n}$. 

Let $\A(\S^{2n-1})$ be the $*$-quotient of $\A(\R^{2n})$ by the two-sided ideal generated by the central element 
$\sum_\mu z^\mu \bar{z}^\mu - 1$. We will denote the images of $z^\mu$ under the quotient map again by $z^\mu$.

A key role in what follows is played by the action of the abelian group $\bT^n$ on $\A(\bR^{2n}_\theta)$ by automorphisms. For $s =(s_\mu) \in \bT^n$, the $*$-automorphism $\sigma_s$ is defined on the generators by $\sigma_s(z^\mu) = e^{2\pi \ii s_\mu} z^\mu$. Clearly, $s \mapsto \sigma_s$ is a group-homomorphism from $\bT^n \to \Aut(\A(\R^{2n}))$. In the special case that $\theta=0$, we see that $\sigma$ is induced by a smooth action of $\bT^n$ on the manifold $\bR^{2n}$. 
Since the ideal generating $\A(\S^{2n-1})$ is invariant under the action of $\bT^n$, $\sigma$ induces a group-homomorphism from $\bT^n$ into the group of automorphisms on the quotient $\A(\S^{2n-1})$ as well. 

We continue by defining the unital $*$-algebra $\A(\R^{2n+1})$ of polynomial functions on the noncommutative $(2n+1)$-plane which is given by adjoining a central self-adjoint generator $x$ to the algebra $\A(\R^{2n})$, {i.e.} $x^*=x$ and $x z^\mu = z^\mu x$ $(\mu=1,\ldots, n)$. The action of the group $\bT^n$ is extended trivially by $\sigma_s(x) = x$. Let $\A(\S^{2n})$ be the $*$-quotient of $\A(\R^{2n+1})$ by the ideal generated by the central element $\sum z^\mu \bar{z}^\mu+x^2 -1$. As before, we will denote the canonical images of $z^\mu$ and $x$ again by $z^\mu$ and $x$, respectively. Since $\bT^n$ leaves this ideal invariant, it induces an action by $*$-automorphisms on the quotient $\A(\S^{2n})$.
\begin{ex}\label{4sp}
For $n=2$ we obtain the noncommutative sphere $S^4_\theta$, which was found in \cite{CL01}. We adopt the notation used therein and let $\A(\S^4)$ be generated by $\alpha,\beta$ and a central $x$ with $x=x^*$ and relations
\be
\alpha \beta = \lambda \beta \alpha, \quad \alpha \beta^* = \bar{\lambda} \beta^* \alpha, \quad \alpha \alpha^* = \alpha^* \alpha, \quad  \beta \beta^* =\beta^* \beta,
\ee
together with the spherical relation $\alpha \alpha^* + \beta \beta^* +x^2 =1$. Here $\lambda=e^{2\pi \ii \theta}$ with $\theta$ a real number. 
\end{ex}
We will now construct a differential calculus on $\R^m$. For $m=2n$, the complex unital associative graded $*$-algebra $\Omega(\R^{2n})$ is generated by $2n$ elements $z^\mu, \bar{z}^\mu$ of degree $0$ and $2n$ elements $dz^\mu,d\bar{z}^\mu$ of degree $1$ with relations:
\bea\nn  \label{rel:diff}
&dz^\mu dz^\nu+ \lambda^{\mu\nu} dz^\nu dz^\mu =0 ;\quad d\bar{z}^\mu dz^\nu+ \lambda^{\nu\mu} dz^\nu d\bar{z}^\mu =0; \quad d\bar{z}^\mu d\bar{z}^\nu+ \lambda^{\mu\nu} d\bar{z}^\nu d\bar{z}^\mu =0;& \\
&z^\mu dz^\nu = \lambda^{\mu\nu} dz^\nu z^\mu;\quad \bar{z}^\mu dz^\nu = \lambda^{\nu\mu} dz^\nu \bar{z}^\mu; \quad \bar{z}^\mu d\bar{z}^\nu = \lambda^{\mu\nu} d\bar{z}^\nu \bar{z}^\mu.&
\eea
There is a unique differential $\dd$ on $\Omega(\R^{2n})$ such that $\dd: z^\mu\mapsto dz^\mu$. The involution $\omega \mapsto \omega^*$ for $\omega \in \Omega(\R^{2n})$ is the graded extension of $z^\mu \mapsto \bar{z}^\mu$, {i.e.} it is such that $(\dd \omega)^*=\dd \omega^*$ and $(\omega_1\omega_2)^* = (-1)^{p_1 p_2}\omega_2^* \omega_1^*$ for $\omega_i\in \Omega^{p_i}(\R^{2n})$.

\noindent
For $m=2n+1$, we adjoin to $\Omega(\R^{2n})$ one generator $x$ of degree $0$ and one generator $dx$ of degree $1$ such that
\be
x dx = dx x; \quad x \omega = \omega x; \quad dx \omega= (-1)^{|\omega|} \omega dx.
\ee
We extend the differential $\dd$ and the graded involution $\omega \mapsto \omega^*$ of $\Omega(\R^{2n})$ to $\Omega(\R^{2n+1})$ by setting $x^*=x$ and $(\dd x)^*=\dd x$, so that $(d x)^*= d x$. 

The differential calculi $\Omega(\S^m)$ on the noncommutative spheres $\S^m$ are defined to be the quotients of $\Omega(\R^{m+1})$ by the differential ideals generated by the central elements $\sum_\mu z^\mu \bar{z}^\mu - 1$ and $\sum z^\mu \bar{z}^\mu+x^2 -1$, for $m=2n-1$ and $m=2n$ respectively. 

The action of $\bT^n$ by $*$-automorphisms on $\A(\M)$ can be easily extended to the differential calculi $\Omega(\M)$, for $M=\R^m$ and $M=\S^m$, by imposing $\sigma_s \circ \dd = \dd \circ \sigma_s$.  

\medskip
In \cite{CD02}, the so-called splitting homomorphism was introduced. For the cases $M=\bR^m$ or $M=S^m$, this map identifies $\A(\M)$ with a subalgebra of $\A(M) \otimes \A(\bT^n_\theta)$, and this identification allows one to use techniques from commutative differential geometry on $\A(M)$ and extend it to $\A(\M)$. Let us recall the definition of the noncommutative $n$-torus $\T^n$. The unital $*$-algebra $\A(\T^n)$ of polynomial functions is generated by $n$ unitary elements $U^\mu$ with relations
\be
U^\mu U^\nu= \lambda^{\mu\nu} U^\nu U^\mu, \quad (\mu,\nu=1,\ldots,n)
\ee
with $\lambda^{\mu\nu}=e^{2\pi\ii\theta_{\mu\nu}}$ as before. There is a natural action of $\bT^n$ on $\A(\T^n)$ by $*$-automorphisms given by $\tau_s (U^\mu)= e^{2 \pi \ii s_\mu } U^\mu$ with $s=(s_\mu) \in \bT^n$. 
This allows one to define a diagonal action $\sigma \times \tau^{-1}$ of $\bT^n$ on $\A(\bR^m\times \T^n):=\A(\bR^m) \otimes \A(\T^n)$ by $s \mapsto \sigma_s \otimes \tau_{-s}$. That is, $s \mapsto (\sigma \times \tau^{-1})_s$ is a group-homomorphisms of $\bT^n$ into $\Aut(\A(\bR^m\times \T^n))$. 

If $\zc^\mu$ denote the classical coordinates of $\bR^n$ corresponding to $z^\mu$ for $\theta=0$, one defines the splitting homomorphism on the generators of $\A(\R^{2n})$ by
\bea
\st: \A(\R^{2n}) &\to& \A(\bR^{2n}) \otimes \A(\T^n); \\ \nn
z^\mu &\mapsto& \zc^\mu \otimes U^\mu.
\eea
One checks that $\st$ induces an isomorphism between the algebra $\A(\R^{2n})$ and the subalgebra $\A(\bR^{2n} \times \T^n)^{\sigma\times\tau^{-1}}$ of $\A(\bR^{2n} \times \T^n)$ consisting of fixed points of the previous diagonal action of $\bT^n$. 
By setting $\st(x)=\xc\otimes 1$ the splitting homomorphism extends trivially to a map from $\A(\R^{2n+1})$ to $\A(\bR^{2n+1}) \otimes \A(\T^n)$, giving an algebra isomorphism $\A(\R^{2n+1})\isom \A(\bR^{2n+1} \times \T^n)^{\sigma\times \tau^{-1}}$.

Furthermore, the map $\st$ will pass to the quotient, for $m=2n,2n+1$,
\be
\st: \A(\S^m) \to \A(S^m) \otimes \A(\T^n)=:\A(S^m \times \T^n),
\ee
giving isomorphisms $\A(\S^m) \isom \A(S^m \times \T^n)^{\sigma\times\tau^{-1}}$. 

The splitting homomorphism allows one to introduce algebras of smooth functions $C^\infty(\M)$, for $M=\bR^m$ or $M=S^m$. They are defined to be the fixed point subalgebras of the diagonal action of $\bT^n$ on $C^\infty(M) \widehat{\otimes} C^\infty(\T^n)$. Here $C^\infty(\T^n)$ is the nuclear Fr\'echet algebra of smooth functions on $\T^n$ and $\widehat{\otimes}$ denotes the completion of the tensor product in the projective tensor product topology (see \cite{CD02} for more details). 

The extension of the splitting homomorphism to the differential calculi yields isomorphisms $\Omega(\M)\isom \big( \Omega(M) \otimes \A(\T^n)\big)^{\sigma \otimes \tau^{-1}}$, with $M$ as above. This allows one to introduce a Hodge star operator on $\Omega(\M)$. Let $\ast$ be the Hodge star operator on $\Omega(M)$ defined with a $\sigma$-invariant metric on $M$. The operator $\ast \otimes \id$ on $\Omega(M) \otimes \A(\T^n)$ restricted to the fixed point subalgebra of the diagonal action, defines the Hodge star operator $\ast_\theta$ on $\Omega(\M)$. Using this operator, one defines a hermitian structure on $\Omega(\M)$ in the following way. If $\omega, \eta \in \Omega^p(\M)$, then
\be \label{def:herm}
\langle \omega, \eta \rangle :=\ast_\theta ( \bar\omega \ast_\theta \eta )
\ee
takes values in $\A(\M)$ and fulfills all properties of a hermitian structure on $\Omega^p(\M)$. 

\medskip
Finally, for the Dirac operator one has the following construction. Suppose for convenience that $M=S^m$ and equip $S^m$ with a Riemannian metric such that $\bT^n$ acts isometrically (this is always possible, for instance by averaging). Let $\cS$ be a spin bundle over the spin manifold $S^m$ and $D$ the Dirac operator on $\Gamma^\infty(S^m,\cS)$. The action of the group $\bT^n$ on $S^m$ does not lift directly to the spinor bundle. Rather, there is a double cover $\pi: \widetilde{\bT}^n \to \bT^n$ and a group-homomorphism $\tilde{s} \to V_{\tilde{s}}$ of $\widetilde{\bT}^n$ into $\Aut(\cS)$ covering the action of $\bT^n$ on $M$:
\be
V_{\tilde{s}} (f \psi) = \sigma_{\pi(s)}(f) V_{\tilde{s}}(\psi),
\ee
for $f \in C^\infty(S^m)$ and $\psi \in \Gamma^\infty(M,\cS)$. It turns out that the proper notion of smooth sections $\Gamma^\infty(\S^m,\cS)$ of a spinor bundle on $\S^m$ is given by the subalgebra of $\Gamma^\infty(S^m,\cS) \widehat\otimes C^\infty(\bT_{\theta/2}^n)$ made of elements which are invariant under the diagonal action 
$V \times \tilde{\tau}^{-1}$ of $\tilde{\bT}^n$. Here $\tilde{s} \mapsto \tilde\tau_{\tilde{s}}$ is the canonical action of $\tilde{\bT}^n$ on $\A(\bT_{\theta/2}^n)$. Since the Dirac operator $D$ will commute with $V_{\tilde{s}}$ one can restrict $D\otimes \id$ to the fixed point subalgebra $\Gamma^\infty(\S^m,\cS)$. 

Next, let $L^2(S^m,\cS)$ be the space of square integrable spinors on $S^m$ and let $L^2(\bT^n_{\theta/2})$ be the completion of $C^\infty(\bT^n_{\theta/2})$ in the norm $f \mapsto \|f\|=\tau(f^*f)^{1/2}$, with $\tau$ the usual trace on $C^\infty(\bT^n_{\theta/2})$. The diagonal action $V \times \tilde{\tau}^{-1}$ of $\tilde{\bT}^n$ extends to $L^2(S^m,\cS)\otimes L^2(\bT^n_{\theta/2})$ and defines $L^2(\S^m, \cS)$ to be the fixed point Hilbert subspace. If $D$ also denotes the closure of the Dirac operator on $L^2(S^m,\cS)$, we denote the operator $D \otimes \id$ on $L^2(S^m,\cS)\otimes L^2(\bT^n_{\theta/2})$ when restricted to $L^2(\S^m, \cS)$ by $D$.

The triple $(C^\infty(\S^m), L^2(\S^m, \cS), D)$ satisfies all axioms of a noncommutative spin geometry (there is also a real structure $J$). In fact, this construction on $\S^m$ can be generalized to any compact Riemannian spin manifold, carrying an isometrical action of $\bT^n$. For more details, we refer to \cite{CL01,CD02}. 

\section{Hopf fibration and associated bundles on $\S^4$} \label{sect:hopf}
We will now construct a $\theta$-deformation of the Hopf fibration $SU(2) \to S^7 \to S^4$. For convenience, the classical fibration is described in some detail in App.~\ref{app:modules}. Firstly, we remind that while there is a $\theta$-deformation of the manifold $S^3\isom SU(2)$, to a sphere $\S^3$, on the latter there is no compatible group structure so that there is no $\theta$-deformation of the group $SU(2)$ \cite{CD02}. Therefore, we must choose the matrix $\theta'_{\mu\nu}$ in such a way that the noncommutative $7$-sphere $\Sk^7$ carries a classical $SU(2)$ action, which in addition is such that the subalgebra of $\A(\Sk^7)$ consisting of $SU(2)$-invariant polynomials is exactly $\A(\S^4)$. As expected, we will find that $\theta'$ is expressed in terms of $\theta$. Then we construct the finitely generated projective modules $\Gamma(\S^4, E^\n)$, associated to the irreducible representations $V^\n$ of $SU(2)$ as the space of $SU(2)$-coequivariant maps from $V^\n$ to $\A(\Sk^7)$.  We will construct projections $p_\n \in \Mat_{4^n}(\A(\S^4))$ such that $\Gamma(\S^4, E^\n) \isom p_\n (\A(\S^4))^{4^n}$. In the special case of the defining representation, we recover the basic instanton projection on the sphere $\S^4$ constructed in \cite{CL01}. 

As mentioned, the interplay of the noncommutative spheres $\S^4$ and $\Sk^7$ is in that $\A(\Sk^7)$ will be required to carry an action of $SU(2)$ by automorphisms and this action is such that
\be \label{inv}
\A(\S^4)= \inv_{SU(2)} (\A(\Sk^7)). 
\ee
These requirements will restrict the values of $\lambda'^{ij}=e^{2\pi\ii \theta'_{ij}}$ is such a manner that there is essentially only `one' noncommutative $7$-sphere such that the invariance condition (\ref{inv}) is satisfied, with a compatible right $SU(2)$ action on $\Sk^7$.
This action 
on the generators of $\A(\Sk^7)$ is simply defined by
\be \label{def:actionSU2}
\alpha_w: (z^1, z^2, z^3, z^4) \mapsto (z^1, z^2, z^3, z^4) 
\begin{pmatrix} w & 0 \\ 0 & w
\end{pmatrix} , \qquad 
w = \begin{pmatrix} w^1 & w^2 \\ -\bar{w}^2 & \bar{w}^1 \end{pmatrix}.
\ee
Here $w^1$ and $w^2$, satisfying $w^1 \bar{w}^1+w^2 \bar{w}^2=1$, are the coordinates on $SU(2)$.
By imposing that the map $w \mapsto \alpha_w$ embeds $SU(2)$ in $\Aut(\A(\Sk^7))$ we find that $\lambda'^{12}=\lambda'^{34} =1$ and $\lambda'^{14}=\lambda'^{23}=\lambda'^{24}=\lambda'^{13}=:\lambda'$. 

In terms of the splitting homomorphism, this means that we can identify $\A(\Sk^7)$ with a certain subalgebra of $\A(S^7 \times \Tk^2)$ instead of $\A(S^7 \times \Tk^4)$. In fact, we can write
\bea 
z^1 = z^1_{(0)} \otimes u, \quad z^3 = z^3_{(0)} \otimes v, \\ \nn
z^2 = z^2_{(0)} \otimes u, \quad z^4 = z^4_{(0)} \otimes v, \nn
\eea
for two unitaries $u, v$ satisfying $uv=\lambda' vu$, i.e. the generators of $\A(\Tk^2)$. 

The subalgebra of $SU(2)$-invariant elements in $\A(\Sk^7)$ can be found in the following way. 
By using the splitting homomorphism, a general element $a\in \A(\Sk^7)$ can be written as a finite sum: $a=\sum a^i_{(0)} \otimes u^i$ where $a^i_{(0)} \in \A(S^7)$ and $u^i \in \A(\Tk^2)$. Then, from the diagonal nature of the action of $SU(2)$ on $\A(\Sk^7)$ and the above formul{\ae} for $z^1, \ldots,z^4$ we have that $\alpha_w(a)= \sum \alpha_w(a^i_{(0)}) \otimes u^i$, encoding the fact that $SU(2)$ essentially acts classically. But this means that any invariant polynomial $a=\alpha_w(a)$ induces a classical invariant polynomial $a_{(0)}$. Hence, the subalgebra of $SU(2)$-invariant elements in $\A(\Sk^7)$ is completely determined by the classical subalgebra of $SU(2)$-invariant elements in $\A(S^7)$. 
From App.~\ref{app:modules} we can conclude  that 
\be
\inv_{SU(2)} (\A(S^7_\theta)) = \bC [~ 1, z^1 \bar{z}^3 +  z^2 \bar{z}^4, -z^1 z^4 + z^2 z^3, z^1 \bar{z}^1+z^2 \bar{z}^2 ~]
\ee
modulo the relations in the algebra $\A(\Sk^7)$. We identify 
\bea \label{subalgebra}
&& \alpha = 2(z^1 \bar{z}^3 +  z^2 \bar{z}^4), \qquad
\beta = 2(-z^1 z^4 + z^2 z^3),  \\ \nn
&& x= z^1 \bar{z}^1+z^2 \bar{z}^2-z^3 \bar{z}^3-z^4 \bar{z}^4,
\eea
and compute that $\alpha \alpha^* + \beta \beta^* + x^2 = 1$. By imposing commutation rules 
$\alpha \beta = \lambda  \beta \alpha$ and $\alpha \beta^* = \bar{\lambda} \beta^* \alpha$, as in Example~\ref{4sp}, we infer that $\lambda'^{14}=\lambda'^{23}=\lambda'^{24}=\lambda'^{13}=\sqrt{\lambda}=:\mu$ on $\Sk^7$. We conclude that $\inv_{SU(2)} (\A(\Sk^7)) = \A(\S^4)$ for $\lambda'^{ij}=e^{2\pi\ii \theta'_{ij}}$ of the following form:
\be
\lambda'_{ij}= \begin{pmatrix} 1 & 1 & \mu & \mu \\ 1 & 1 & \mu & \mu \\
 \bar{\mu} & \bar{\mu} &1 & 1\\ \bar{\mu} & \bar{\mu} &1 & 1 \end{pmatrix}, \quad \mu = \sqrt{\lambda}, 
\ee
or equivalently 
\be
\theta'_{ij}=\frac{\theta}{2}\begin{pmatrix} 0 & 0 & 1 & 1 \\ 0 & 0 & 1 & 1 \\ -1 & -1 & 0 & 0 \\  -1 & -1 & 0 & 0\end{pmatrix}.
\ee	

There is a nice description of the instanton projection constructed in \cite{CL01} in terms of ket-valued  polynomials on $\Sk^7$. The latter are elements in the right $\A(\Sk^7)$-module $\cE:=\bC^4 \otimes \A(\Sk^7)=: \A(\Sk^7)^4$ with a hermitian structure given by $\langle \xi, \eta\rangle=\sum_j \xi^*_j \eta_j$. To any $\ket{\xi}\in \cE$ one associates its dual $\bra{\xi}\in \cE^*$ by setting $\bra{\xi}(\eta):= \langle \xi,\eta \rangle$, $\forall \eta\in\cE$. 

Similarly to the classical case (see App.~\ref{app:modules}), we define $\ket{\psi_1},\ket{\psi_2} \in \A(\Sk^7)^4$ by
\be \label{def:kets}
\ket{\psi_1}=(z^1 , -\bar{z}^2, z^3, -\bar{z}^4 )^\t , \quad 
\ket{\psi_2}=(z^2 , \bar{z}^1, z^4, \bar{z}^3 )^\t , 
\ee
with $\t$ denoting transposition. They satisfy $\langle\psi_k |\psi_l \rangle = \delta_{kl}$, so that the $4\times 4$-matrix $p=\ket{\psi_1} \bra{\psi_1} + \ket{\psi_2} \bra{\psi_2}$ is a projection, $p^2=p=p^*$, with entries in $\A(\S^4)$. Indeed, let us introduce the matrix 
\be 
u = (\ket{\psi_1},\ket{\psi_2}) =
\begin{pmatrix} 
z^1 & z^2 \\ 
-\bar{z}^2 & \bar{z}^1 \\
z^3 & z^4 \\
-\bar{z}^4 & \bar{z}^3
\end{pmatrix}.
\ee
Then $u^* u = \I_2 $ and $p = u u^*$.
The action \eqref{def:actionSU2} becomes 
\be \label{actionSU2}
\alpha_w (u) = u w,
\ee
from which the invariance of the entries of $p$ follows at once. Explicitly one finds 
\be \label{projection1}
p=\frac{1}{2}\begin{pmatrix}
1+x & 0 & \alpha & \beta \\
0 & 1+x &-\mu \beta^* & \bar{\mu} \alpha^* \\
\alpha^*& -\bar{\mu} \beta & 1-x & 0\\
\beta^* & \mu \alpha & 0 & 1-x 
\end{pmatrix}.
\ee
The projection $p$ is easily seen to be equivalent to the projection describing the instanton on $\S^4$ constructed in \cite{CL01}. Indeed, if one defines 
\be  \label{def:ketsequiv}
\ket{\widetilde{\psi}_1}=(z^1 , -\mu \bar{z}^2, z^3, -\bar{z}^4 )^\t , \quad 
\ket{\widetilde{\psi}_2}=(z^2 , \mu \bar{z}^1, z^4, \bar{z}^3 )^\t , 
\ee
one obtains exactly the projection obtained therein, that is, 
\be
\widetilde{p}=\frac{1}{2}\begin{pmatrix}
1+x & 0 & \alpha & \beta \\
0 & 1+x &-\lambda \beta^* & \alpha^* \\
\alpha^*& -\bar{\lambda} \beta & 1-x & 0\\
\beta^* & \alpha & 0 & 1-x 
\end{pmatrix}
\ee

\medskip
We will denote the image of $p$ in $\A(\S^4)^4$ by $\Gamma(\S^4, E)=p \A(\S^4)^4$ which is clearly a right $\A(\S^4)$-module.
Another description of the module $\Gamma(\S^4,E)$ comes from considering coequivariant maps from $\bC^2$ to $\A(\Sk^7)$ \cite{Dur96}. The defining left representation of $SU(2)$ on $\bC^2$ is given by $SU(2) \times \bC^2 \to \bC^2; (w,v) \mapsto w \cdot v$. The collection $\Hom_{SU(2)} (\bC^2, \A(\Sk^7))$ of coequivariant maps, {i.e.} of maps $\phi: \bC^2 \to \A(\Sk^7)$, such that
\be \label{def:equiv}
\phi(w^{-1} \cdot v) = \alpha_w (\phi(v)), 
\ee
is a right $\A(\S^4)$-module (it is in fact also a left $\A(\S^4)$-module).

Since $SU(2)$ acts classically on $\A(\Sk^7)$, one sees that the coequivariant maps are given on the canonical basis $\{e_1,e_2\}$ of $\bC^2$ by $\phi(e_k) = \langle\psi_k|f\rangle$ for $\ket{f} = \ket{f_1,f_2,f_3,f_4}^\t$, with $f_i \in \A(\S^4)$ (cf. App.~\ref{app:modules}). We then have the following isomorphism
\bea
\Gamma(\S^4, E) &\isom& \Hom_{SU(2)} (\bC^2, \A(\Sk^7)) \\ \nn
\sigma= p \ket{f} &\leftrightarrow& \phi: e_k \mapsto\langle\psi_k|f\rangle.
\eea

More generally, one can define the right $\A(\S^4)$-module $\Gamma(\S^4, E^\n)$ associated with any irreducible representation $\rho_n : SU(2) \to GL(V^\n)$, with $V^\n = \Sym^n (\bC^2)$, for a positive integer $n$. The module of coequivariant maps $\Hom_{\rho_n} (V^\n, \A(\Sk^7))$ consists of maps $\phi:V^\n \to \A(\Sk^7)$ satisfying
\be \label{def:equivgeneral}
\phi(\rho_n^{-1} (w) \cdot v) = \alpha_w (\phi(v)).
\ee
It is easy to see that these maps are of the form $\phi_\n(e_k) = \langle\phi_k|f\rangle$ on the basis $\{e_1, \ldots, e_{n+1}\}$ of $V^\n$ where now $\ket{f}\in \A(\S^4)^{4^n}$ and 
\be \label{defphi}
\ket{\phi_k} = \frac{1}{a_k} \ket{\psi_1}^{\otimes(n-k+1)} \otimes_S \ket{\psi_2}^{\otimes(k-1)}  \quad(k=1,\ldots,n+1),
\ee
with $\otimes_S$ denoting symmetrization  and $a_k$ are suitable normalization constants. These vectors $\ket{\phi_k} \in \bC^{4^n} \otimes \A(\Sk^7)=: \A(\Sk^7)^{4^n}$ are orthogonal (with the natural hermitian structure), and with $a_k^2=\binom{n}{k-1}$ they are also normalized. Then 
\be \label{def:proj}
p_\n:=\ket{\phi_1} \bra{\phi_1} + \ket{\phi_2} \bra{\phi_2}+\cdots + \ket{\phi_{n+1}} \bra{\phi_{n+1}} \in \Mat_{4^n}(\A(\S^4))
\ee
defines a projection $p^2=p=p^*$. That its entries are in $\A(\S^4)$ and not in $\A(\Sk^7)$ is easily seen. Indeed, much as it happens for the vector $u$ in \eqref{actionSU2}, for every $i=1,\ldots,4^n$, the vector $u_{(i)}=\big( \ket{\phi_1}_i ,  \ket{\phi_2}_i, \ldots , \ket{\phi_{n+1}}_i \big)$ transforms under the action of $SU(2)$ to the vector $\big( \ket{\phi_1}_i , \ldots , \ket{\phi_{n+1}}_i \big) \cdot \rho_{(n)} (w)$ so that each entry $\sum_k \ket{\phi_k}_i \bra{\phi_k}_j$ of $p_\n$ is $SU(2)$-invariant and hence an element in $\A(\S^4)$. With this we proved the following.
\hyphenation{iso-mor-phic}
\begin{prop} \label{prop:modulesisomorphism}
The module of coequivariant maps $\Hom_{\rho_n} (V^\n, \A(\Sk^7))$ is isomorphic to $\Gamma(\S^4,E^\n):=p_\n (\A(\S^4)^{4^n})$ (as right-$\A(\S^4)$ modules) with the isomorphism given explicitly by:
\begin{eqnarray*}
\Gamma(\S^4,E^\n) &\isom& \Hom_{\rho_n} (V^\n, \A(\Sk^7)) \\ \nn
\sigma_\n= p_\n \ket{f} &\leftrightarrow& \phi_\n: e_k \mapsto \langle\phi_k|f\rangle.
\end{eqnarray*}
\end{prop}
Using the splitting homomorphisms of the previous section, one can lift this whole construction to the smooth level. One proves that the $C^\infty(\S^4)$-module $\Gamma^\infty(\S^4,E^\n)$ defined by $p_\n (C^\infty(\S^4))^{4^n}$ is isomorphic to $\Hom_{\rho_n} (V^\n, C^\infty(\Sk^7))$.

With the projections $p_\n$ one associates (Grassmannian) connections on the modules $\Gamma(\S^4, E^\n)$ in a canonical way:
\be \label{def:gras}
\nabla = p_\n \circ \dd : \Gamma(\S^4, E^\n) \to \Gamma(\S^4, E^\n) \otimes_{\A(\S^4)} \Omega^1(\S^4)
\ee
where $(\Omega^*(\S^4), \dd)$ is the differential calculus defined in the previous section. 
An expression for these connections as acting on coequivariant maps can be obtained using the above isomorphism and results in:
\be \label{def:grasconn}
\nabla (\phi)(e_k) =\dd (\phi(e_k)) + A_{kl} \phi(e_l)
\ee
where $A_{kl}=\langle \phi_k | \dd \phi_l \rangle \in \Omega^1(\Sk^7)$. The corresponding matrix 
$A$ is called the connection one-form; it is clearly anti-hermitian, and it is valued in the derived representation space, $\rho'_n : su(2) \to \End(V^\n)$, of the Lie algebra $su(2)$. 

\medskip
The case $n=1$ describes classically ({i.e.} $\theta =0$) the charge $-1$ instanton \cite{Ati79}. An instanton is defined as a connection on $\Gamma(S^4,E)$ with (anti-)selfdual curvature $F$, {i.e} $\ast F = \pm F$ with $\ast$ the Hodge star operator. In physics, instantons are of importance since they are extrema of the Yang-Mills action. The equation of motion obtained from this action by a variatonal method is called the Yang-Mills equation $[\nabla,\ast F]=0$. In the case that $F$ is (anti-)selfdual, this equation of motion follows directly from the Bianchi identity $[\nabla, F]=0$. There is no need to stress the huge importance of instantons (and in general Yang-Mills gauge theory) both in physics and in mathematics.

On the noncommutative sphere $\S^4$, the curvature of the connection $p \circ \dd$ constructed above satisfies the following anti-selfdual equation \cite{CD02}
(see also \cite{AB02,Lnd01})\footnote{An early attempt to write self-duality equations in terms of projections was done in \cite{Dub83}.}:
\be
\ast_\theta p (\dd p)^2 = -p (\dd p)^2.
\ee

In order to fully justify the name instanton, one should find a noncommutative analogue of the Yang-Mills action such that connections with an (anti-)selfdual curvature are its extrema. This will be discussed elsewhere \cite{LS04}. 

\section{Index of Dirac operators} \label{sect:index}
We know from Sect.~\ref{sect:NCS} that there is a structure of noncommutative spin geometry
on the sphere $\S^4$ given by a `triple' $(C^\infty(\S^4), L^2(\S^4, \cS), D, \gamma)$ with $\gamma=\gamma_5$ the grading. In this Section we shall compute explicitly the index of the Dirac operator with coefficients in the bundles $E^\n$, that is the index of the operator 
of $D_{p_\n}:= p_\n (D \otimes \I_{4^n})p_\n$. In order to do that, we will use the (`even dimensional' version of the) local index formula of Connes and Moscovici \cite{CM95} which we shall briefly describe. 

Suppose in general that $(\A,\cH,D,\gamma)$ is an even $p$-summable spectral triple with discrete simple dimension spectrum. 
Let $C_*(\A)$ be the complex consisting of cycles over the algebra $\A$, that is in degree $n$, $C_n(\A):=\A^{\otimes(n+1)}$. On this complex there are defined the Hochschild operator $b:C_n(\A)\to C_{n-1}(\A)$ and the boundary operator $B:C_n(\A)\to C_{n+1}(\A)$, satisfying $b^2=0, B^2=0, bB+Bb=0$; thus $(b+B)^2=0$. From general homological theory, one defines a bicomplex $CC_*(\A)$ by $CC_{(n,m)}(\A):= CC_{n-m}(\A)$ in bi-degree $(n,m)$. Dually, one defines $CC^\ast(\A)$ as functionals on $CC_*(\A)$, equipped with the dual Hochschild operator $b$ and coboundary operator $B$ (we refer to \cite{C94} and \cite{Lod92} for more details on this). 
\begin{thm}[Connes-Moscovici \cite{CM95}]~ \\
\begin{itemize}
\item[(a)] An even cocycle $\phi^*=\sum_{k\geq0} \phi^k$ in $CC^*(\A)$, $(b+B)\phi^*=0$, defined by the following formul{\ae}. For $k=0$,
\be
\phi^0(a):=\resz  ~z^{-1} \tr( \gamma a |D|^{-2z} );
\ee
whereas for $k\neq 0$
\be\label{evencocy}
\phi^{2k}(a^0, \ldots, a^{2k}) := \sum_{\alpha} c_{k,\alpha} \resz ~ \tr \big( \gamma a^0 [D,a^1]^{(\alpha_1)} \cdots [D,a^{2k}]^{(\alpha_{2k})} |D|^{-2(|\alpha| +k + z)} \big)
\ee
where 
\bd
c_{k,\alpha}=(-1)^{|\alpha|} \Gamma(k+|\alpha|) \big( \alpha! (\alpha_1+1) (\alpha_1+\alpha_2+2) \cdots (\alpha_1+ \cdots +\alpha_{2k} + 2k)\big)^{-1}
\ed 
and $T^{(j)}$ denotes the j'th iteration of the derivation $T \mapsto [D^2,T]$. 
\item[(b)]
For $e \in K_0(\A)$, the Chern character $\ch_*(e) = \sum_{k\geq 0} \ch_k(e)$ is the even cycle in $CC_*(\A)$, $(b+B)\ch_*(e)=0$, defined by the following formul{\ae}. For $k=0$,
\be
\ch_0(e):= \tr(e);
\ee
whereas for $k\neq 0$
\be \label{chernchar}
\ch_k(e):=(-1)^k \frac{(2k)!}{k!}  \sum (e_{i_0 i_1}-\frac{1}{2}\delta_{i_0 i_1}) \otimes {e_{i_1 i_2} \otimes e_{i_1 i_2} \otimes \cdots \otimes e_{i_{2k} i_0} }.
\ee
\item[(c)] The index is given by the natural pairing between cycles and cocycles
\be
\ind D_{e}=\langle \phi^\ast,\ch_\ast(e) \rangle .
\ee
\end{itemize}
\end{thm}

\medskip
We concentrate on a compact Riemannian spin manifold $M$ of even dimension carrying an isometric action of an $n$-torus. Set $\cH:=L^2(M,\cS)$ and recall the grading on $\cB(\cH)$ with respect to the action of $\bT^n$ \cite{CL01}. An element $T\in \cB(\cH)$ that is smooth for the action of $\bT^n$, ({i.e.} such that the map $\bT^n \ni s \mapsto \alpha_s(T)$ with $\alpha_s$ defined by $\alpha_s(T) := U(s) T U(s)^{-1}$ is smooth for the norm topology) can be expanded as $T =\sum T_{r}$ with $r=(r_1,r_2,\ldots,r_n)$ a multi-index, and with each $T_{r}$ of homogeneous degree  $r$ under the action of $\bT^n$, {i.e.}
\bd
\alpha_s(T_{r}) =e^{2 \pi \ii(\sum_{\mu=1}^n r_\mu s_\mu )} T_{r}  \quad (s \in \bT^n) .
\ed
Note that $\alpha_s$ coincides on $\pi(C^\infty(M)) \subset \cB(\cH)$ with the automorphism $\sigma_s$ defined in Sect.~\ref{sect:NCS}.
Then, let ($p_1, p_2,\ldots, p_n$) be the infinitesimal generators of the action of $\bT^n$ so that $U(s)=\exp{2 \pi \ii(\sum_{\mu=1}^n s_\mu p_\mu)}$.  For $T \in \cB(\cH)$ we define a twisted representation on $\cH$ by
\be\label{twist}
L_\theta(T):=\sum_r T_r U(r_\mu \theta_{\mu 1}, \ldots, r_\mu \theta_{\mu n} ) = \sum_r T_r \exp \big\{ 2 \pi \ii \sum_\mu r_\mu \theta_{\mu\nu} p_\nu \big\}
\ee
with $\theta$ an $n\times n$ anti-symmetric matrix. Since $\bT^n$ acts by isometries, each $p_\mu$ commutes with $D$ so that the latter is of degree $0$ and $L_\theta([D,a])=[D,L_\theta(a)]$ for $a \in C^\infty(M)$. It was shown in \cite{CL01} (to which we refer for more details) that $(L_\theta(C^\infty(M)), \cH, D)$ satisfies  all axioms of a noncommutative spin geometry (there is also a grading $\gamma=\gamma_5$ and a real structure $J$). In fact, the algebra $L_\theta(C^\infty(M))$ is isomorphic to the algebra $C^\infty(M_\theta)$ of the previous section. 

\medskip
As a next step, we write the $\phi^{2k}$ that define the local index formula in \eqref{evencocy}, 
by means of the twist $L_\theta$. Let $f^0,\ldots, f^{2k}\in C^\infty(M)$ and suppose that the operator 
$f^0 [D,f^1] \cdots [D,f^{2k}]$ is a homogeneous element of degree $r$. Then, as in \eqref{twist} 
\be
L_\theta(f^0 [D,f^1] \cdots [D,f^{2k}]) = f^0 [D,f^1] \cdots [D,f^{2k}] U(r_\mu \theta_{\mu1},\ldots, r_\mu \theta_{\mu n}) .
\ee
Each term in the local index formula for $(L_\theta(C^\infty(M)),\cH,D)$ then takes the form
\bd
\resz  \tr \big( \gamma f^0 [D,f^1]^{(\alpha_1)} \cdots [D,f^{2k}]^{(\alpha_{2k})} |D|^{-2(|\alpha| +k + z)} U(s) \big)
\ed
for $s_\nu=r_\mu \theta_{\mu\nu}$ so that $s \in \bT^n$. The appearance of $U(s)$ here, is a consequence of the close relation with the index formula for a $\bT^n$-equivariant Dirac spectral triple on $M$. In \cite{CH97}, Chern and Hu considered an even dimensional compact spin manifold $M$ on which a (connected compact) Lie group $G$ acts by isometries. The equivariant Chern character was defined as an equivariant version of the JLO-cocycle, the latter being an element in equivariant entire cyclic cohomology. The essential point is that they obtained an explicit formula for the above residues. In the case of the previous $\bT^n$-action on $M$, one gets
\begin{multline}
\resz  \tr \big( \gamma f^0 [D,f^1]^{(\alpha_1)} \cdots [D,f^{2k}]^{(\alpha_{2k})} |D|^{-2(|\alpha| +k + z)} U(s) \big) \\
=\Gamma(|\alpha|+k) \lim_{t \to 0} t^{|\alpha|+k}  \tr \big( \gamma f^0 [D,f^1]^{(\alpha_1)} \cdots [D,f^{2k}]^{(\alpha_{2k})} e^{-tD^2} U(s) \big)
\end{multline}
for every $s \in \bT^n$; moreover, this limit vanishes when $|\alpha|\neq0$ (Thm 2 in \cite{CH97}). Combining these results, we arrive at the following Lemma. 
\begin{lma}
Let $(L_\theta(C^\infty(M)),\cH,D)$ be the spectral triple defined above. Then all terms in $\phi^*$ with $|\alpha|\neq 0$ vanish and the local index formula takes the form:
\be
\phi^{2k}(a^0, \ldots, a^{2k}) = c_k \resz  \tr \big( \gamma a^0 [D,a^1] \cdots [D,a^{2k}] |D|^{-2(k + z)} \big)
\ee
where $c_k=(k-1)! / (2k)!$. \end{lma}

In our case of interest, the index of the Dirac operator on $\S^4$ with coefficients in some noncommutative vector bundle determined by $e \in K_0(C(\S^4))$, we obtain
\bea
\ind D_{e}=\langle\phi^*,\ch_*(e)\rangle &=&\resz z^{-1} \tr \big( \gamma \pi_D (\ch_0(e)) |D|^{-2z} \big) \nn \\
&&+\frac{1}{2!} \resz  \tr \big( \gamma \pi_D (\ch_1(e)) |D|^{-2-2z} \big) \nn  \\
&&+\frac{1}{4!} \resz  \tr \big( \gamma \pi_D (\ch_2(e)) |D|^{-4-2z} \big)
\eea
Here $\pi_D$ is the representation of the universal differential calculus given by
\be
\pi_D : \Omega^p_\un(\A(\S^4)) \to \cB(\cH) , \quad 
a^0 \delta a^1 \cdots \delta a^p \mapsto a^0 [D,a^1] \cdots [D,a^p].
\ee
Let us examine at which quotients of $\Omega_\un(\A(\S^4))$ this representation $\pi_D$ is well-defined. Unfortunately, $\pi_D$ is not well-defined on the quotient $\Omega(\S^4)$ defined in the previous section. For example already $[D,\alpha][D,\alpha] \neq 0$ whereas $d\alpha d\alpha=0$ in $\Omega(\S^4)$. This was already noted in \cite{CD02} and in fact 
\be
\Omega(\S^4) \isom \pi_D \big(\Omega_\un(\A(\S^4)) \big) / {\pi_D(\delta J_0)}
\ee
where $J_0:= \{ \omega \in \Omega_\un(\A(\S^4)) | \pi_D(\omega)=0 \}$ are the so-called 'junk-forms' \cite{C94}. We will avoid a discussion on junk-forms and introduce instead a different quotient of $\Omega_\un(\A(\S^4))$. We define $\Omega_D(\S^4)$ to be $\Omega_\un(\A(\S^4))$ modulo the relations 
\begin{align}
&\alpha \delta \beta - \lambda (\delta \beta) \alpha = 0, \quad
(\delta \alpha) \beta - \lambda \beta \delta \alpha = 0, \nn \\
&\alpha \delta \beta^* - \bar \lambda (\delta \beta^*) \alpha = 0, \quad 
(\delta \alpha^*) \beta - \bar \lambda \beta \delta \alpha^* = 0, \nn \\
& a \delta x - (\delta x) a = 0 , \qquad \forall a \in \A(\S^4),
\end{align}
avoiding the second order relations that define $\Omega(\S^4)$. Using the splitting homomorphism one proves that the above relations are in the kernel of $\pi_D$, for instance, 
$\alpha [D,\beta]-\lambda [D,\beta] \alpha = 0$ so that $\pi_D$ is well-defined on $\Omega_D(\S^4)$.

In App.~\ref{app:chern} we compute the Chern characters as elements in $\Omega_D(\S^4)$, which results in the following Lemma. 
\begin{lma}
The following formul{\ae} hold for the images under $\pi_D$ of the Chern characters of $p_\n$:
\begin{eqnarray*}
\pi_D(\ch_0(p_\n))&=& n+1;\\
\pi_D(\ch_1(p_\n))&=&0;\\
\pi_D(\ch_2(p_\n))&=&\frac{1}{6} n(n+1)(n+2)\pi_D( \ch_2(p_{(1)})); \\
\end{eqnarray*}
up to the coefficients $\mu_k=(-1)^k \frac{(2k)!}{k!}$.
\hfill $\Box$
\end{lma}
Combining this with the simple form of the index formula while taking the proper coefficients, we find that
\be
\ind D_{p_\n} = \frac{1}{4!} \frac{4!}{2!} \frac{1}{6} n(n+1)(n+2) \resz  \tr\big(\gamma \pi_D(\ch_2(p_{(1)})) |D|^{-4-2z} \big)
\ee
where for the vanishing of the first term, we used the fact that $\ind D=0$, since the first Pontrjagin class on $S^4$ vanishes. Thm I.2 in \cite{CM95} allows one to express the residue as a Dixmier trace. Combining this with $\pi_D(\ch_2(p_{(1)}))=3 \gamma$ (as computed in \cite{CL01}), we obtain
\bd
3 \cdot \resz  \tr (|D|^{-4-2z} ) = 6 \cdot \dix( |D|^{-4} )= 2
\ed
since the Dixmier trace of $|D|^{-m}$ on the $m$-sphere equals $8/m!$ (cf. for instance \cite{GVF01,Lnd97}). This combines to give:
\begin{prop}
The index of the Dirac operator on $\S^4$ with coefficients in $E^\n$ is given by:
\bd
\ind D_{p_\n} = \frac{1}{6} n(n+1)(n+2).
\ed
\rightbox
\end{prop}
\noindent Note that this coincides with the classical result.
\section{The noncommutative principal bundle} \label{sect:hopf-galois}
In this section, we apply the general theory of Hopf-Galois extensions \cite{KT81,Mon93} to the inclusion $\A(\S^4) \into \A(\Sk^7)$. Such extensions can be understood as noncommutative principal bundles. We will first dualize the construction of the previous section, {i.e.} replace the action of $SU(2)$ on $\A(\Sk^7)$ by a coaction of $\A(SU(2))$. Then, we will recall some definitions involving Hopf-Galois extensions and principality (\cite{BH04}) of such extensions. We show that $\A(\S^4) \into \A(\Sk^7)$ is a not-cleft (i.e. not-trivial) principal Hopf-Galois extension and compare the connections on the associated bundles, induced from the strong connection, with the Grassmannian connection defined in Sect.~\ref{sect:hopf}.

The action of $SU(2)$ on $\A(\Sk^7)$ by automorphisms can be easily dualized to a coaction $\Delta_R: \A(\Sk^7) \to \A(\Sk^7) \otimes \A(SU(2))$, where now $\A(SU(2))$ is the unital complex $*$-algebra generated by $w^1, \bar{w}^1,w^2, \bar{w}^2$ with relation $w^1 \bar{w}^1+w^2 \bar{w}^2=1$. Clearly, $\A(SU(2))$ is a Hopf algebra with comultiplication
\be
\Delta : \begin{pmatrix} w^1 &w^2\\-\bar{w}^2 & \bar{w}^1 \end{pmatrix} \mapsto  \begin{pmatrix} w^1 &w^2\\-\bar{w}^2 & \bar{w}^1 \end{pmatrix} \otimes  \begin{pmatrix} w^1 &w^2\\-\bar{w}^2 & \bar{w}^1 \end{pmatrix},
\ee
antipode $S(w^1) = \bar{w}^1, S(w^2)= -w^2$ and counit $\epsilon(w^1)=\epsilon(\bar{w}^1)=1, \epsilon(w^2)=\epsilon(\bar{w}^2)=0$. The coaction of $\A(SU(2))$ on $\A(\Sk^7)$ is given by
\be \label{def:coactionSU2}
\Delta_R: (z^1, z^2, z^3, z^4) \mapsto (z^1, z^2, z^3, z^4) \otimes \begin{pmatrix} w^1 & w^2 & 0 & 0 \\ -\bar{w}^2 & \bar{w}^1  & 0 & 0 \\ 0 & 0 & w^1 & w^2 \\ 0 & 0 & -\bar{w}^2 & \bar{w}^1 \end{pmatrix}.
\ee
The algebra of coinvariants in $\A(\Sk^7)$, which consists of elements $p\in \A(\Sk^7)$ satisfying $\Delta_R(p)= p \otimes 1$, can be identified with $\A(\S^4)$ for the particular values of $\theta'_{ij}$ found before, in the same way as in Sect.~\ref{sect:hopf}.

The associated modules $\Gamma(\S^4,E^\n)$ are described in the following way. Given an irreducible corepresentation of $\A(SU(2))$, $\rho_\n : V^\n \to \A(SU(2))\otimes V^\n$  with $V^\n = \Sym^n (\bC^2)$, we  denote $\rho_\n(v)=v_{(0)} \otimes v_{(1)}$. Then, the module of coequivariant maps $\Hom^{\rho_\n} (V^\n, \A(\Sk^7))$ consists of maps $\phi:V^\n \to \A(\Sk^7)$ satisfying
\be \label{def:coequiv}
\phi(v_{(1)}) \otimes S v_{(0)}=\Delta_R \phi(v); \quad v \in \bC^2.
\ee
Again, such maps are $\bC$-linear maps of the form $\phi_\n(e_k) = \langle\phi_k|f\rangle$ on the basis $\{e_1, \ldots, e_{n+1}\}$ of $V^\n$ in the notation of the previous section. Also, Proposition \ref{prop:modulesisomorphism} above translates straightforwardly into the isomorphism $\Hom^{\rho_\n} (V^\n, \A(\Sk^7)) \isom p_\n (\A(\S^4))^{4^n}$ for the projections defined in equation (\ref{def:proj}). 

Before we proceed, recall that for an algebra $P$ and a subalgebra $B \subset P$, $P\otimes_B P$ denotes the quotient of the tensor product $P \otimes P$ by the ideal generated by expressions $p\otimes b p'- p b \otimes p'$, for $p,p'\in P, b \in B$. 
\begin{defn}
Let $H$ be a Hopf algebra and $P$ a right $H$-comodule algebra, {i.e.} such that the coaction $\Delta_R : P \to P \otimes H$ is an algebra map. Let the algebra of coinvariants be $B:=\coinv_{\Delta_R}(P):=\{p\in P:\Delta_R(p)=p\otimes 1$ \}. One says that $B \into P$ is a \textbf{Hopf-Galois extension} if the canonical map 
\be\label{can}
\chi : P \otb P \to P \otimes H; \quad p' \otimes_B p \mapsto p' \Delta_R(p) = p' p_{(0)} \otimes p_{(1)}
\ee
is bijective. 
\end{defn}
We use Sweedler-like  notation for the coaction: $\Delta_R(p)= p_{(0)} \otimes p_{(1)}$.
The canonical map is left $P$-linear and right $H$-colinear and is a morphism (an isomorphism  for Hopf-Galois extensions) of left $P$-modules and right $H$-comodules. It is also clear that $P$ is 
both a left and a right $B$-module.

Classically, the notion of Hopf-Galois extension corresponds to freeness of the action of a Lie group $G$ on a manifold $P$. Indeed, freeness can be translated into bijectivity of the map
\be
\tilde{\chi}: P \times G \to P \times_G P ,  \quad
(p,g) \mapsto (p,p \cdot g), 
\ee 
where $P\times_G P$ denotes the fibred direct product consisting of elements $(p,p')$ with the same image under the quotient map $P \to P/G$. 

\medskip
For a Hopf algebra $H$ which is cosemisimple, surjectivity of the canonical map \eqref{can} implies its bijectivity \cite{Sch90}. Moreover, in order to prove surjectivity of $\chi$, it is enough to prove that for any generator $h$ of $H$, the element $1\otimes h$ is in the image of the canonical map. Indeed, if 
$\chi(g_k \otb g'_k) =1 \otimes g$ and $\chi(h_l \otb h'_l) =1 \otimes h$ for $g,h\in H$, 
then $\chi(g_k h_l \otb h_l' g_k')=g_k h_l \chi(1 \otb h_l' g_k')= 1\otimes h g$, using the fact that the canonical map restricted to $1\otimes_B P$ is a homomorphism. Extension to all of $P\otb P$ then follows from left $P$-linearity of $\chi$. It would also be easy to write down an explicit expression for the inverse of the canonical map. Indeed, one has $\chi^{-1} (1 \otimes hg) = g_k h_l \otimes_B h_l' g_k'$ in the above notation so that the general form of the inverse follows again from left $P$-linearity. 

\begin{prop} \label{prop:hopfgalois}
The inclusion $\A(\S^4) \into \A(\Sk^7)$ is a Hopf-Galois extension. 
\end{prop}
\noindent
\textit{Proof. }Since $\A(SU(2))$ is cosemisimple, we can rely for a proof of this statement on 
the previous remarks. On the other hand, it is straightforward to check that in terms of the ket-valued polynomials defined in \eqref{def:kets} we have
\bean
&\chi(\sum_i \bra{\psi_1}_i \ot4 \ket{\psi_1}_i)=1\otimes w^1; 
& \quad \chi(\sum_i \bra{\psi_1}_i \ot4 \ket{\psi_2}_i)=1\otimes w^2; \\ 
&\chi(\sum_i \bra{\psi_2}_i \ot4 \ket{\psi_1}_i)=-1\otimes \bar{w}^2;  & \quad \chi(\sum_i \bra{\psi_2}_i \ot4 \ket{\psi_2}_i)=1\otimes \bar{w}^1.
\eean
\rightbox

In the definition of a principal bundle in differential geometry there is much more than the requirement of bijectivity of the canonical map. It turns out that our `structure group' being $H=\A(SU(2))$ which, besides being cosemisimple has also bijective antipode, all additional desired properties follows from the surjectivity of the canonical map
which we have just established.
We refer to \cite{SS04,BH04} for the full fledged theory while giving only the basic definitions that we shall need.

For our purposes, a better algebraic translation of the notion of a principal bundle is  encoded in the 
requirement that the extension $B\subset P$, besides being Hopf-Galois,  is also 
faithfully flat. We recall \cite{lam} that 
a right module $P$ over a ring $R$ is  said to be \textbf{faithfully flat} 
if the functor $P \otimes_R \cdot$  is exact and faithful  on the category 
$_R \mathcal{M}$ of left $R$-modules.
\noindent Flatness means that the functor associates exact sequences of abelian
groups to exact sequences of $R$-modules and the functor is faithful if it is
injective on morphisms.
Equivalently one could state that a right module $P$ over a ring $R$ is
faithfully flat if a sequence $M' \to M \to M''$ in $_R \mathcal{M}$ is exact if and only if  $P 
\otimes_R M' \to P \otimes_R M \to P \otimes_R M''$ is exact.

As mentioned, from the fact that $H=\A(SU(2))$ is both cosemisimple and has also bijective antipode, 
the faithful flatness of  $\A(\Sk^7)$ as a right (as well as left) $\A(\S^4)$-module follows from the surjectivity of the canonical map (\cite{Sch90}, Th. I). \\

One says that a principal Hopf-Galois extension is \textbf{cleft} if there exists a (unital) convolution-invertible colinear map $\phi: H \to P$, called a \textbf{cleaving map} \cite{DGH01,SS04}. Classically, this notion is close (although not equivalent) to triviality of a principal bundle \cite{DHS99}. In \cite{BM93} (cf. \cite{HM99}) it is shown that if a principal Hopf-Galois extension is cleft, its associated modules are trivial, {i.e.} isomorphic to the free module $B^N$ for some $N$. In our case, we can conclude the following. 
\begin{prop}
The Hopf-Galois extension $\A(\S^4) \into \A(\Sk^7)$ is not cleft.
\end{prop}
\noindent
{\it Proof.}
This is a simple consequence of the nontriviality of the Chern characters of the projection $p_\n$ as seen in Sect.~\ref{sect:index}. Indeed, this implies that the associated modules are nontrivial. $\hfill \Box$

Summing up what we have shown up to now, we have the following
\begin{thm}
The inclusion $\A(\S^4) \into \A(\Sk^7)$ is a not-cleft faithfully flat $\A(SU(2))$-Hopf-Galois extension. 
\end{thm}

An important consequence is the existence of a so-called \textbf{strong connection} \cite{Haj96,DGH01}. In fact, the existence of such a connection could be used to give a more intuitive definition of `principality of an extension' \cite{BH04}.
Let us first recall that if $H$ is cosemisimple and has a bijective antipode, then a $H$-Hopf-Galois extension $B \into P$ is \textbf{equivariantly projective}, that is, there exists a left $B$-linear right $H$-colinear splitting $s:P\to B\otimes P$ of the multiplication map $m:B \otimes P \to P$, $m\circ s=\id_P$ \cite{SS04}. Such a map characterizes a strong connection. 

\begin{defn}
Let $B \into P$ be a $H$-Hopf-Galois extension. A \textbf{strong connection one-form} is a 
map $\omega: H \to \Omega^1_\un P$ satisfying
\begin{enumerate}
\item $\bar\chi \circ \omega = 1 \otimes (\id -\epsilon)$, \qquad  (\textit{fundamental vector field condition})
\item $\Delta_{\Omega^1_\un (P)} \circ \omega = (\omega \otimes \id) \circ \Ad_R$, \qquad 
(\textit{right adjoint colinearity})
\item $\delta p - p_{(0)} \omega(p_{(1)}) \in ( \Omega^1_\un B) P$, \, $\forall  p \in P$, 
\qquad (\textit{strongness condition}).
\end{enumerate}
\end{defn}
Here $\Delta_R : P \to P \otimes H$, $\Delta_R(p)= p_{(0)} \otimes p_{(1)}$, is extended to  $\Delta_{\Omega^1_\un (P)}$ on $\Omega^1_\un P \subset P \otimes P$ in a natural way by 
\be
\Delta_{\Omega^1_\un (P)} (p'\otimes p) \mapsto p'_{(0)} \otimes p_{(0)} \otimes p'_{(1)} p_{(1)}, 
\ee
and $\Ad_R(h) = h_{(2)}\otimes S(h_{(1)}) h_{(3)}$ is the right adjoint coaction of $H$. Finally, the map $\bar\chi: P\otimes P \to P \otimes H$ is defined like the canonical map as $\bar\chi(p'\otimes p)=p'p_{(0)}\otimes p_{(1)}$. 

As shown in \cite{BH04} (cf. \cite{BDZ04,HMS01}), a strong connection can always be given by a map 
$\ell: H \to P \otimes P$ satisfying
\bea \label{ell}
\ell(1) &=& 1\otimes 1 , \nn \\ 
\bar\chi(\ell(h))&=&1\otimes h ,  \nn \\
(\ell \otimes \id)\circ \Delta &=& (\id \otimes \Delta_R) \circ \ell , \nn \\ 
(\id \otimes \ell)\circ \Delta &=& (\Delta_L \otimes \id) \circ \ell ,
\eea
where $\Delta_L:P \to H \otimes P$, $p \mapsto S^{-1} p_{(1)} \otimes p_{(0)}$. Then, one defines the 
connection one-form by
\be
\omega:h \mapsto \ell(h) - \epsilon(h) 1\otimes 1.
\ee
Indeed, if one writes $\ell(h)=h^{\br{1}} \otimes h^{\br{2}}$ (summation understood) and applies $\id \otimes \epsilon$ to the second formula in \eqref{ell}, one has $h^{\br{1}}h^{\br{2}}=\epsilon(h)$. Therefore, 
\be
\omega( h)= h^{\br{1}} \delta  h^{\br{2}}
\ee
where $\delta: P \to \Omega^1_\un P, ~p \mapsto 1\otimes p -p \otimes 1$. 
Equivariant projectivity of $B\into P$ follows by taking as splitting of the multiplication the map
$s:P \to B\otimes P,  ~p \to p_{(0)} \ell(p_{(1)})$. 

For later use, we prove the following Lemma, analogous to the strongness condition 3. above. 
\begin{lma} \label{leftstrongness}
Let $\omega$ be a strong connection one-form on a $H$-Hopf-Galois extension $B \into P$ with the antipode of $H$ invertible. Then
$$
\delta p + \omega(S^{-1}p_{(1)}) p_{(0)} \in P \Omega_\un^1 B, \quad \forall p \in P.
$$
\end{lma}
\noindent 
\textit{Proof.} By writing $\omega$ in terms of $\ell$ it follows that $\delta p + \omega(S^{-1}p_{(1)}) p_{(0)}$ reduces to the expression $-p\otimes 1 + l(S^{-1} p_{(1)}) p_{(0)}$. From the second property of $\ell$ in \eqref{ell}, it follows that this expression is in the kernel of $\bar\chi$. Since $\chi$ is an isomorphism, $\delta p + \omega(S^{-1}p_{(1)}) p_{(0)}$ is in the ideal generated by expressions of the form $p \otimes b p' - p b \otimes p'$. In other words, it is an element in $P \Omega_\un^1(B) P$. Finally, it is not difficult to show that
$$
(\id \otimes \Delta_R ) \big( \delta p + \omega(S^{-1}p_{(1)}) p_{(0)} \big) = \big( \delta p + \omega(S^{-1}p_{(1)}) p_{(0)}\big) \otimes 1
$$
from which we conclude that $\delta p + \omega(S^{-1}p_{(1)}) p_{(0)}$ is in fact in $P \Omega_\un^1(B)$. $\hfill \Box$

\medskip
In our case, the existence of a strong connection follows from \cite{SS04}. However, we will write an explicit expression in terms of the inverse of the canonical map. If we denote the latter when lifted to $P \otimes P$ by $\tau$ it follows that $\ell(h) = \tau(1\otimes h)$ satisfies the same recursive relation found before for $\chi^{-1}$ (proof of Proposition \ref{prop:hopfgalois} above): if $\ell(h)=h_l \otimes h_l'$ and $\ell(g)=g_k \otimes g_k'$, then 
\be \label{recursive}
\ell(hg)=g_k h_l \otimes h_l' g_k' .
\ee
It turns out that in our case the map $\ell:H \to P\otimes P$ defined in this way defines a strong connection. 
\begin{prop}
On the Hopf-Galois extension $\A(\S^4) \into \A(\Sk^7)$, the following formul{\ae} on the generators of 
$\A(SU(2))$, 
\bea \label{def:lGenerators}
&\ell(w^1) = \sum_i \bra{\psi_1}_i \otimes \ket{\psi_1}_i; \quad \ell(w^2)= \sum_i \bra{\psi_1}_i \otimes \ket{\psi_2}_i;& \\ \nn
&\ell(\bar{w}^2)=-\sum_i \bra{\psi_2}_i \otimes \ket{\psi_1}_i; \quad \ell(\bar{w}^1)=\sum_i \bra{\psi_2}_i \otimes \ket{\psi_2}_i.&  \nn
\eea
define a strong connection.
\end{prop}
\noindent
{\it Proof.} We extend the expressions \eqref{def:lGenerators} to all of $\A(SU(2))$ by giving recursive relations, using formula (\ref{recursive}). Recall the usual vector basis $\{ r^{klm} : k \in \bZ, m,n \geq 0\}$ in $\A(SU(2))$ given by
\be
r^{klm}:=\left\{ \begin{array}{ll} (-1)^n (w^1)^k(w^2)^m(\bar{w}^2)^n & k \geq 0,\\
(-1)^n (w^2)^m(\bar{w}^2)^n (\bar{w}^1)^{-k} & k<0. \end{array} \right.
\ee
The recursive expressions on this basis are explicitly given by
\bea \nn \label{recursiveEll}
\ell(r^{k+1,mn}) &=& \bar{z}^1 \ell(r^{kmn}) z^1 +z^2 \ell(r^{kmn}) \bar{z}^2+\bar{z}^3 \ell(r^{kmn}) z^3+z^4 \ell(r^{kmn}) \bar{z}^4,  \quad k\geq 0, \\ \nonumber
\ell(w^{k-1,mn}) &=&  \bar{z}^2 \ell(r^{kmn}) z^2+z^1 \ell(r^{kmn}) \bar{z}^1 +\bar{z}^4 \ell(r^{kmn}) z^4+ z^3 \ell(r^{kmn}) \bar{z}^3,  \quad k <0, \\ \nonumber
\ell(w^{k,m+1,n}) &=&  \bar{z}^1 \ell(r^{kmn}) z^2-z^2 \ell(r^{kmn}) \bar{z}^1 +\bar{z}^3 \ell(r^{kmn}) z^4 -z^4 \ell(r^{kmn}) \bar{z}^3, \\
\ell(w^{km,n+1}) &=& \, \bar{z}^2 \ell(r^{kmn}) z^1 -z^1 \ell(r^{kmn}) \bar{z}^2+\bar{z}^4 \ell(r^{kmn}) z^3-z^3 \ell(r^{kmn}) \bar{z}^4, 
\eea
while setting $\ell(1)=1 \otimes 1$. In essentially the same manner as was done in \cite{BCDT03} (although much simpler in our case) we prove that $\ell$ defined by the above recursive relations indeed satisfies all conditions of a strong connection. \hfill $\Box$

\medskip
The strong connection on the extension $\A(\S^4) \into \A(\Sk^7)$ induces connections on the associated modules in the following way \cite{HM99}. For $\phi \in \Hom^{\rho_\n} (V^\n, \A(\Sk^7))$, we set 
\be
\nabla_\omega (\phi)(v) \mapsto \delta \phi(v) + \omega(v_{(0)}) \phi(v_{(1)}).
\ee
Using the right adjoint colinearity of $\omega$ and a little algebra one shows that $\nabla_\omega (\phi)$ satisfies the following coequivariance condition
$$
\nabla_\omega (\phi)(v_{(1)}) \otimes S v_{(0)} = \Delta_{\Omega^1_\un(P)} \big( \nabla_\omega (\phi)(v) \big)
$$
so that 
$$\nabla_\omega: \Hom^{\rho_\n} (V^\n, \A(\Sk^7)) \to \Hom^{\rho_\n} (V^\n, \Omega_\un^1(\A(\Sk^7))).$$ 
In fact, from Lemma \ref{leftstrongness} it follows that $\nabla_\omega$ is a map from $\Hom^{\rho_\n} (V^\n, \A(\Sk^7))$ to $\Hom^{\rho_\n} (V^\n, \A(\Sk^7)) \otimes \Omega^1_\un(\A(\S^4))$. This allows one to compare it to the Grassmannian connection of equation (\ref{def:gras}). It turns out that the connection one-form $\omega$ coincides with the connection one-form $A$ of equation (\ref{def:grasconn}), on the quotient $\Omega^1(\Sk^7)$ of $\Omega^1_\un(\A(\Sk^7))$. More precisely, let $\{e^\n_k\}$ be a basis of $V^\n$, and $e^\n_{kl}$ the corresponding matrix coefficients of $\A(SU(2))$ in the representation $\rho_\n$. An explicit expression for $\omega(e^{(n)}_{kl})$ can be obtained from equations \eqref{recursiveEll}; for example $\omega(e^{(1)}_{kl})=\braket{\psi_k}{\delta \psi_l}, \, k,l=1,2$.

By using these and formul{\ae} \eqref{expr:gauge1}-\eqref{gauge2}, one shows that 
$$\pi(\omega(e^\n_{kl}))=A^\n_{kl}= \braket{\phi_k}{\dd \phi_l},$$ 
where $\pi: \Omega_\un (\A(\Sk^7)) \to \Omega(\Sk^7)$ is the quotient map.

\appendix
\section{Associated modules}\label{app:modules}
We will review the classical construction of the instanton bundle on $S^4$ \cite{Ati79} taking the approach of \cite{Lnd00}. We generalize slightly and construct complex vector bundles on $S^4$ associated to all finite-dimensional irreducible representations of $SU(2)$. 

We start by recalling the Hopf fibration $\pi:S^7 \to S^4$. Let 
\bea
S^7&:=&\{z=(z^1,z^2, z^3,z^4):|z^1|^2+|z^2|^2+|z^3|^2+|z^4|^2=1\}, \nn \\
S^4&:=&\{(\alpha,\beta,x):\alpha \bar{\alpha} + \beta \bar{\beta} +x^2=1\}, \nn \\
SU(2)&:=&\{ w \in GL(2,\bC) : w^* w=w w^* =1, \det w=1\} \nn \\ 
&=& \left\{ w=\begin{pmatrix} w^1 & w^2 \\-\bar{w}^2 & \bar{w}^1\end{pmatrix} :
w^1 \bar{w}^1+w^2 \bar{w}^2=1 \right\}.
\eea
The space $S^7$ carries a right $SU(2)$-action:
\be
S^7 \times SU(2) \to S^7 , \quad 
\big(z,w\big) \mapsto (z^1, z^2, z^3, z^4) \begin{pmatrix} w  & 0 \\ 0 &w \end{pmatrix}.
\ee
The Hopf map is defined as a map $\pi(z^1, z^2, z^3, z^4)\mapsto( \alpha, \beta,x)$ where
\bea \label{hopfmap}
&& \alpha = 2(z^1 \bar{z}^3 +  z^2 \bar{z}^4), \quad
\beta = 2(-z^1 z^4 + z^2 z^3),\\ \nn
&& x= z^1 \bar{z}^1+z^2 \bar{z}^2-z^3 \bar{z}^3-z^4 \bar{z}^4 , \nn
\eea
and one computes $\alpha \bar{\alpha} + \beta \bar{\beta} +x^2=(\sum_i|z^i|^2)^2=1$. 

The finite-dimensional irreducible representations of $SU(2)$ are labeled by a positive integer $n$ with 
$n+1$-dimensional representation space $V^\n \isom \Sym^n(\bC^2)$. The space of smooth $SU(2)$-equivariant maps from $S^7$ to $V^\n$ is defined by
\be
C_{SU(2)}^\infty(S^7, V^\n):=\{\phi: S^7 \to V^\n : \phi(z \cdot w) = w^{-1} \cdot \phi(z)\}.
\ee
We will now construct projections $p_\n$ as $N\times N$ matrices taking values in $C^\infty(S^4))$, such that $\Gamma^\infty(S^4, E^\n):=p_\n C^\infty(S^4)^N $ is isomorphic to $C_{SU(2)}^\infty(S^7, V^\n)$ as right $C^\infty(S^4)$-modules. As the notation suggests, $E^\n$ is the vector bundle over $S^4$ associated with the corresponding representation. Let us first recall the case $n=1$ from \cite{Lnd00} and then use this to generate the vector bundles for any $n$. The $SU(2)$-equivariant maps from $S^7$ to $V^{(1)} \isom \bC^2$ are of the form
\be \label{equivariantmaps}
\phi_{(1)}(z) = \begin{pmatrix} \bar{z}^1 \\ \bar{z}^2\end{pmatrix} f_1 +\begin{pmatrix} -z^2 \\ z^1\end{pmatrix} f_2 +\begin{pmatrix} \bar{z}^3 \\ \bar{z}^4\end{pmatrix} f_3+ \begin{pmatrix} -z^4 \\ z^3\end{pmatrix} f_4,
\ee
where $f_1, \ldots,f_4$ are smooth functions that are invariant under the action of $SU(2)$, {i.e.} they are functions on the base space $S^4$. 

A nice description of the equivariant maps is given in terms of ket-valued functions $\ket{\xi}$ on $S^7$, which are then elements in the free module $\cE:=\bC^N \otimes C^\infty(S^7)$. The $C^\infty(S^7)$-valued hermitian structure on $\cE$ given by $\langle \xi, \eta\rangle=\sum_j \xi^*_j \eta_j$ allows one to associate dual elements $\bra{\xi}\in \cE^*$ to each $\ket{\xi}\in \cE$ by $\bra{\xi}(\eta):= \langle \xi,\eta \rangle, ~\forall \eta\in\cE$. 

If we define $\ket{\psi_1},\ket{\psi_2} \in \A(S^7)^4$ by
\be 
\ket{\psi_1}=(z^1 , -\bar{z}^2, z^3, -\bar{z}^4 )^\t;\quad 
\ket{\psi_2}=(z^2 , \bar{z}^1, z^4, \bar{z}^3 )^\t,
\ee
with $\t$ denoting transposition, the equivariant maps in (\ref{equivariantmaps}) are given by 
\be
\phi_{(1)}(z) = \begin{pmatrix} \langle\psi_1|f\rangle \\ \langle\psi_2|f\rangle \end{pmatrix}, 
\ee
where $\ket{f} \in (C^\infty(S^4))^4 := \bC^4 \otimes C^\infty(S^4)$. 
Since $\langle \psi_k | \psi_l \rangle =\delta_{kl}$ as is easily seen, we can define a projection in $M_4(C^\infty(S^4))$ by 
\be
p_{(1)}=\ket{\psi_1} \bra{\psi_1} + \ket{\psi_2} \bra{\psi_2}.
\ee
Indeed, by explicit computation we find a matrix with entries in $C^\infty(S^4)$ which is the limit of the projection (\ref{projection1}) for $\theta=0$. Denoting the right $C^\infty(S^4)$-module $p_{(1)} (C^\infty(S^4))^4$ by $\Gamma(S^4, E^{(1)})$, we have
\bea
\Gamma(S^4, E^{(1)}) &\isom& C_{SU(2)}^\infty(S^7, \bC^2) \\ \nn
\sigma_{(1)}= p_{(1)} \ket{f} &\leftrightarrow& \phi_{(1)}=\begin{pmatrix} \langle\psi_1|f\rangle \\ \langle\psi_2|f\rangle \end{pmatrix}.
\eea

For the general case, we note that the $SU(2)$-equivariant maps from $S^7$ to $V^\n$ are of the form
\be
\phi_\n(z) = \begin{pmatrix} \langle\phi_1|f\rangle \\ \vdots \\ \langle\phi_{l+1}|f\rangle \end{pmatrix}
\ee
where $\ket{f}\in C^\infty(S^4)^{4^n}$ and $\ket{\phi_k}$ is the completely symmetrized form of the tensor product $\ket{\psi_1}^{\otimes n-k+1} \otimes \ket{\psi_2}^{\otimes k-1}$ for $k=1,\ldots,n+1$, normalized to have norm 1 as in formula (\ref{defphi}). For example, for the adjoint representation $n=2$, we have
\bea
\ket{\phi_1}&:=& \ket{\psi_1} \otimes \ket{\psi_1},\\ \nn
\ket{\phi_2}&:=& \frac{1}{\sqrt{2}} \big(\ket{\psi_1} \otimes \ket{\psi_2}+\ket{\psi_2} \otimes \ket{\psi_1}\big),\\ \nn
\ket{\phi_3}&:=& \ket{\psi_2} \otimes \ket{\psi_2}.
\eea
Since in general, $\langle \phi_k | \phi_l\rangle=\delta_{kl}$, the matrix-valued function
$$
p_\n=\ket{\phi_1} \bra{\phi_1} + \ket{\phi_2} \bra{\phi_2}+\cdots  +\ket{\phi_{n+1}} \bra{\phi_{n+1}} \in M_{4^n}(C^\infty(S^4))
$$
defines a projection whose entries are in $C^\infty(S^4)$, since each entry $\sum_k \ket{\phi_k}_i \bra{\phi_k}_j$ is $SU(2)$-invariant (cf. below formula (\ref{def:proj})). We conclude that
\bea
p_\n (C^\infty(S^4)^{4^n}) &\isom& C_{SU(2)}^\infty(S^7, V^\n) \\ \nn
\sigma_\n= p_\n \ket{f} &\leftrightarrow& \phi_\n=\begin{pmatrix} \langle\phi_1|f\rangle \\ \vdots \\ \langle\phi_{n+1}|f\rangle \end{pmatrix}.
\eea

\section{Chern characters} \label{app:chern}
We will compute the Chern characters of the projections $p_\n$ defined in Sect.~\ref{sect:hopf}. It turns out to be sufficient for our purposes to obtain expressions in the differential calculus $\Omega_D(\S^4) \subset \Omega_D(\Sk^7)$ defined in Sect.~\ref{sect:hopf}, which is a quotient of the universal differential calculus.

In the definition of the Chern character (\ref{chernchar}) we can replace the tensor product by the universal differential $\delta$ by the isomorphism
$$
\A \otimes \bar{\A}^{\otimes p} \isom \Omega^p_\un(\A),
$$
where $\bar{\A}:=\A/{\bC \I}$ and $\Omega_\un(\A) = \oplus_p \Omega^p_\un(\A)$ is the universal differential algebra generated by 
$a \in \A$  and symbols $\delta a$, $a\in \A$ of order 1 satisfying 
$$\delta(ab)=(\delta a) b+ a \delta b \quad \delta (\alpha a +\beta b) = \alpha \delta a + \beta \delta b; \quad (a,b\in\A,\alpha,\beta\in \bC).
$$
Then 
\be
\ch_k(e):= (-1)^k \frac{(2k)!}{k!} \big\langle (e-\frac{1}{2}) (\delta e)^{2k} \big\rangle \in \Omega^{2k}_\un(\A).
\ee
We recall the differential calculi $\Omega_D(\S^4)$ and $\Omega_D(\Sk^7)$ from Sect.~\ref{sect:hopf}. We defined $\Omega_D(\S^4)$ as the quotient of $\Omega_\un(\A(\S^4))$ by the relations
\begin{align}
&\alpha \delta \beta =\lambda (\delta \beta) \alpha, \quad
(\delta \alpha) \beta = \lambda \beta \delta \alpha, \nn \\
&\alpha \delta \beta^* = \bar \lambda (\delta \beta^*) \alpha, \quad 
(\delta \alpha^*) \beta = \bar \lambda \beta \delta \alpha^*, \nn \\
&a \delta x = (\delta x) a , \qquad (a \in \A(\S^4)).
\end{align}
The inclusion $\A(\S^4) \into \A(\Sk^7)$ extends to an injective map $\Omega_D(\S^4) \to \Omega_D(\Sk^7)$, where $\Omega_D(\Sk^7)$ is the quotient of $\Omega_\un(\A(\Sk^7))$ by only the relations in (\ref{rel:diff}) of order one, that is by the relations:
\be
z^i \delta z^j= \lambda^{ij} (\delta z^j)z^i; \quad z^i \delta \bar{z}^j= \lambda^{ji} (\delta \bar{z}^j) z^i .
\ee
Recall that the projections $p_\n$ were defined by $p_\n = \sum_k \ket{\phi_k} \bra{\phi_k}$, 
where $\ket{\phi_k}$ with $k=1,\ldots,n+1,$ is given by
\be
\ket{\phi_k} = \frac{1}{a_k} \ket{\psi_1}^{\otimes(n-k+1)} \otimes_S \ket{\psi_2}^{\otimes(k-1)}, 
\quad a_k^2=\binom{n}{k-1}.
\ee

Before we start the computation of the Chern characters, we state the computation rules in $\Omega_D(\Sk^7)$. Firstly, from the very definition of the vectors $\ket{\phi_k}$ and the inner product in $\cE \otimes_\bC \cE \otimes_\bC \cdots \otimes_\bC \cE$, we can express, for any $k=1,\ldots, n+1$,
\bea \label{expr:gauge1}
\langle \phi_k | \delta \phi_{k-1} \rangle &=& \sqrt{(n-k)(k+1)} \langle \psi_2 | \delta \psi_1 \rangle \\ 
\langle \phi_k | \delta \phi_{k+1} \rangle &=& \sqrt{(n-k-1)(k+2)} \langle \psi_1 | \delta \psi_2 \rangle 
 \label{expr:gauge3} \\
\langle \phi_k | \delta \phi_k \rangle &=& (n-k-1) \langle \psi_1 | \delta \psi_1 \rangle + (k+1) 
\langle \psi_2 | \delta \psi_2 \rangle \nn \\ &=& (n-2k-2) \langle \psi_1 | \delta \psi_1 \rangle ,\label{gauge2}
\eea
by using the relation $\langle \psi_2 | \delta \psi_2 \rangle = - \langle \psi_1 | \delta \psi_1 \rangle .$  The previous are in fact the only nonzero expressions for $\braket{\phi_k}{\delta \phi_l}$. If we apply $\delta$ to these equations, we obtain expressions for $\langle \delta \phi_k | \delta \phi_l \rangle$ in terms of $\psi_1$ and $\psi_2$. From this, we deduce the following result that will be central in the computation of the Chern characters. 
\begin{lma} \label{lemma:techn}
The following relations hold in $\Omega_D(\Sk^7)$:
\bean 
\sum_{k,l=1}^{n+1} \braket{\phi_k}{\delta \phi_l} \braket{\phi_l}{\delta \phi_k}&=& \frac{1}{6}n (n+1)(n+2) \sum_{r,s=1}^2 \braket{\psi_r}{\delta \psi_s} \braket{\psi_s}{\delta \psi_r}, \\
~ \\
\sum_{k,l,m=1}^{n+1} \braket{\phi_k}{\delta \phi_l} \braket{\phi_l}{\delta \phi_m}  \braket{\phi_m}{\delta \phi_k}
&=& \frac{1}{6}n (n+1)(n+2) \sum_{r,s,t=1}^2 \braket{\psi_r}{\delta \psi_s} \braket{\psi_s}{\delta \psi_t}  \braket{\psi_t}{\delta \psi_r}.
\eean
\end{lma}
\noindent Of course, there will be similar formul{\ae} for $\braket{\delta \phi_k}{\delta \phi_l} \braket{\phi_l}{\delta \phi_k}$, etc. 

\medskip
The zeroeth Chern character is easy to compute: 
\be
\ch_0(p_\n) = \tr (p_\n) = \sum_k \langle \phi_k | \phi_k \rangle= n+1 .
\ee
In the computation of $\ch_1(p_\n)$ we use the relation $\langle \delta \phi_k | \phi_l \rangle = - \langle \phi_k | \delta \phi_l \rangle$, which follows from applying the derivation $\delta$ to $\langle \phi_k | \phi_l \rangle = \delta_{kl}$ and the fact that in $\Omega_D(\Sk^7)$, $\braket{\phi_k}{\delta \phi_l}$ commutes with any element in $\A(\Sk^7)$, in particular with $\bra{\phi_m}_i$. Thus, 
\bea \nn
&&\ch_1(p) = \sum \ket{\phi_k} \braket{\phi_k}{\delta \phi_l} \braket{\phi_l}{\delta \phi_m} \bra{\phi_m} + 
\sum \ket{\phi_k} \ketbra{\delta \phi_k}{\delta\phi_m} \bra{\phi_m} \\ \nn
&& \qquad - \frac{1}{2} \sum \bigg( \ket{\delta \phi_k} \braket{\phi_k}{\delta \phi_l} \bra{\phi_l}+\ket{\delta \phi_k} \bra{\delta \phi_k} + \ket{\phi_k} \braket{\delta\phi_k}{\delta \phi_l} \bra{\phi_l} + \ket{\phi_k}\braket{\delta \phi_k}{\phi_l} \bra{\delta \phi_l} \bigg) \\ \nn
&& \qquad  = \frac{1}{2} 
\sum_{m=1}^{n+1}\big( \braket{\delta \phi_m}{\delta \phi_m} - \ket{\delta \phi_m}\bra{\delta\phi_m} \big)
\eea
 By using equation (\ref{gauge2}) and its analogue for $\ket{\delta \phi_m}\bra{\delta \phi_m}$, $m=1,\ldots,n+1$,
\bd
\ket{\delta \phi_m}\bra{\delta \phi_m}=(k+1) \langle \psi_1 | \delta \psi_1 \rangle + (n-k-1) \langle \psi_2 | \delta \psi_2 \rangle,
\ed
we find that
\be
\ch_1(p_\n) = \frac{1}{2} n(n+1) \big( \langle \psi_1 | \delta \psi_1 \rangle + \langle \psi_2 | \delta \psi_2 \rangle \big)=\frac{1}{2}n(n+1) \ch_1(p_{(1)}).
\ee
Note that this equation holds in the differential subalgebra $\Omega_D(\S^4)$. Since $\ch_1(p_{(1)})$ was shown to vanish in \cite{CL01}, we proved the vanishing of the first Chern character in $\Omega_D(\S^4)$. The vanishing of $\ch_1(p_{(1)})$ can also be seen from the explicit form of $\ket{\psi_1}$ and $\ket{\psi_2}$.

A slightly more involved computation in $\Omega_D(\Sk^7)$ shows that
\bea \nn
\ch_2(p_\n) &=& \frac{1}{2} \sum \big\{ \delta \big(\braket{\phi_k}{\delta \phi_l}\braket{\phi_l}{\delta \phi_m}\braket{\phi_m}{\delta \phi_k} \big ) + \braket{\delta \phi_k}{\delta \phi_l}\braket{\phi_l}{\delta \phi_m}\braket{\phi_m}{\delta \phi_k} \\
&& + \braket{\delta \phi_k}{\delta \phi_l}\braket{\delta \phi_l}{\delta \phi_k} + \delta \big( \braket{\delta \phi_k}{\delta \phi_l}\braket{\phi_l}{\delta \phi_k} \big) \big\}.
\eea
And by using Lemma~\ref{lemma:techn} we finally get
\be
\ch_2(p_\n) = \frac{1}{6} n (n+1)(n+2) \ch_2(p_{(1)}), 
\ee
as an element in $\Omega_D^4(\S^4)$.

\subsection*{Acknowledgments} 
We are grateful to Tomasz Brzezi{\'n}ski, Ludwik D\c{a}browski, Piotr M. Hajac,
Cesare Reina, Chiara Pagani for very useful discussions.


\begin{thebibliography}{90}

\bibitem{AB02}
P.~Aschieri and F.~Bonechi.
\newblock On the noncommutative geometry of twisted spheres.
\newblock Lett. Math. Phys. 59 (2002) 133-156.

\bibitem{Ati79}
M.~F. Atiyah.
\newblock {\em The Geometry of Yang-Mills Fields}.
\newblock Fermi Lectures. Scuola Normale Pisa, 1979.

\bibitem{ADHM78}
M.~F. Atiyah, N.~J. Hitchin, V.~G. Drinfeld and Yu.~I. Manin.
\newblock Construction of instantons.
\newblock Phys. Lett. A65 (1978) 185-187.

\bibitem{BCDT03}
F.~Bonechi, N.~Ciccoli, L.~D\c{a}browski and M.~Tarlini.
\newblock Bijectivity of the canonical map for the noncommutative instanton
  bundle.
\newblock J. Geom. Phys. 51 (2004) 419-432.

\bibitem{BDZ04}
T.~Brzezi\'nski, L.~D\c{a}browski and B.~Zielinski.
\newblock {H}opf fibration and monopole connection over the contact quantum
  spheres.
\newblock J. Geom. Phys. 51 (2004) 71-81.

\bibitem{BH04}
T.~Brzezi\'nski and P.~M. Hajac.
\newblock The {C}hern-{G}alois character.
\newblock C.R. Acad. Sci. Paris Ser. I 338 (2004) 113-116.

\bibitem{BM93}
T.~Brzezi\'nski and S.~Majid.
\newblock Quantum group gauge theory on quantum spaces.
\newblock Commun. Math. Phys. 157 (1993) 591-638.

\bibitem{CH97}
S.~Chern and X.~Hu.
\newblock Equivariant {C}hern character for the invariant {D}irac operator.
\newblock Michigan Math. J. 44 (1997) 451-473.

\bibitem{C94}
A.~Connes.
\newblock {\em Noncommutative Geometry}.
\newblock Academic Press, San Diego, 1994.

\bibitem{CD02}
A.~Connes and M.~Dubois-Violette.
\newblock Noncommutative finite-dimensional manifolds. {I}. {S}pherical
  manifolds and related examples.
\newblock Commun. Math. Phys. 230 (2002) 539-579.

\bibitem{CL01}
A.~Connes and G.~Landi.
\newblock Noncommutative manifolds: {T}he instanton algebra and isospectral
  deformations.
\newblock Commun. Math. Phys. 221 (2001) 141-159.

\bibitem{CM95}
A.~Connes and H.~Moscovici.
\newblock The local index formula in noncommutative geometry.
\newblock Geom. Funct. Anal. 5 (1995) 174-243.

\bibitem{DGH01}
L.~D\c{a}browski, H.~Grosse and P.~M. Hajac.
\newblock Strong connections and {C}hern-{C}onnes pairing in the
  {H}opf-{G}alois theory.
\newblock Commun. Math. Phys. 220 (2001) 301-331.

\bibitem{DHS99}
L.~D\c{a}browski, P.~M. Hajac and P.~Siniscalco.
\newblock Explicit {H}opf-{G}alois description of {$SL_{e^{2\pi
  i/3}}(2)$}--induced {F}robenius homomorphisms.
\newblock In D.~Kastler et.~al., editors, {\em Enlarged
  Proceedings of the ISI GUCCIA Workshop on Quantum Groups, Noncommutative
  Geometry and Fundamental Physical Interactions}. Nova Science Pub., Commack-NY, 
  1999, pp. 279-298.

\bibitem{Dub83}
M. Dubois-Violette, {\it Equations de Yang et Mills, mod\`eles $\sigma$ \`a deux dimensions et g\'en\'eralisation}, in `{\it Math\'ematique et Physique}', Progress in Mathematics, vol. 37, Birkh\"auser 1983, pp.~43-64.\\M. Dubois-Violette and Y. Georgelin, {\it Gauge Theory in Terms of Projector Valued Fields}, Phys. Lett. 82B (1979) 251-254. 

\bibitem{Dur96}
M.~Durdevich.
\newblock Geometry of quantum principal bundles {I}.
\newblock Commun. Math. Phys. 175 (1996) 427-521.\\
M.~Durdevich.
\newblock Geometry of quantum principal bundles {II}.
\newblock Rev. Math. Phys. 9 (1997) 531-607.

\bibitem{GVF01}
J.~M. Gracia-Bond{\'\i}a, J.~C. V\'arilly and H.~Figueroa.
\newblock {\em Elements of {N}oncommutative {G}eometry}.
\newblock Birkh\"auser, Boston, 2001.

\bibitem{Haj96}
P.~M. Hajac.
\newblock Strong connections on quantum principal bundles.
\newblock Commun. Math. Phys. 182 (1996) 579--617.

\bibitem{HM99}
P.~M. Hajac and S.~Majid.
\newblock Projective module description of the {$q$}-monopole.
\newblock Commun. Math. Phys. 206 (1999) 247-264.

\bibitem{HMS01}
P.~M. Hajac, R.~Matthes and W.~Szyma\'nski.
\newblock A locally trivial quantum {H}opf fibration.
\newblock \texttt{arXiv:math.QA/0112317}, to appear in {A}lgebra and
  {R}epresentation {T}heory.

\bibitem{KT81}
H.~F. Kreimer and M.~Takeuchi.
\newblock {H}opf algebras and {G}alois extensions of an algebra.
\newblock Indiana Univ. Math. J. 30 (1981) 675-692.

\bibitem{lam}
T.~Y. Lam.
\newblock {\em Lectures on modules and rings}.
\newblock Spring-Verlag, New-York, 1999

\bibitem{Lnd97}
G.~Landi.
\newblock {\em An Introduction to Noncommutative
  Spaces and their Geometry}.
\newblock Springer-Verlag, Berlin 1997.

\bibitem{Lnd00}
G.~Landi.
\newblock Deconstructing monopoles and instantons.
\newblock Rev. Math. Phys. 12 (2000) 1367-1390.

\bibitem{Lnd01} G. Landi, talk at the Mini-workshop on {\it Noncommutative Geometry BetweenMathematics and Physics}, Ancona, February 23-24, 2001. 

\bibitem{LS04}
G.~Landi and W.~van Suijlekom.
\newblock In preparation.

\bibitem{Lod92}
J.-L. Loday.
\newblock {\em Cyclic Homology}.
  Springer-Verlag, Berlin, 1992.

\bibitem{Mon93}
S.~Montgomery.
\newblock {\em Hopf algebras and their actions on rings}.
\newblock AMS, 1993.

\bibitem{Rie90a}
M.~A. Rieffel.
\newblock Non-commutative tori - {A} case study of non-commutative
  differentiable manifolds.
\newblock Contemp. Math. 105 (1990) 191-212.

\bibitem{SS04}
P.~Schauenburg and H.-J. Schneider.
\newblock Galois type extensions and {H}opf algebras.
\newblock To be published.

\bibitem{Sch90}
H.-J. Schneider.
\newblock Principal homogeneous spaces for arbitrary {H}opf algebras.
\newblock Israel J. Math. 72 (1990) 167-195.

\end{thebibliography}

\end{document}